\documentclass[english,aps,a4paper,twocolumn,showpacs,showkeys,amsfonts,amsmath,amssymb,amscd,nofootinbib]{revtex4}
\usepackage[english]{babel}
\usepackage{graphicx}
\usepackage[dvips]{color}

\def\Journal#1#2#3#4#5#6{{#1}, {\it #4} \textbf{#5}, #6 (#2).}
\def\Book#1#2#3#4#5{{#1}, {\it #3} (#4, #5, #2).}

\newcommand{\zrp}{z_+}
\newcommand{\zrm}{z_-}

\newcommand{\beq}{\begin{equation}}
\newcommand{\enq}{\end{equation}}
\newcommand{\es}{\end{subequations}}

\newcommand{\la}{\lambda}

\newcommand{\EE}{\mathbb{E}}
\newcommand{\TT}{\mathbb{T}}\newcommand{\FF}{\mathbb{F}}
\newcommand{\RR}{\mathbb{R}}

\def\eq{\equiv}

\def\fo{\forall}

\def\p{\partial}

\def\un{\underset}

\begin{document}

\title{Exit times in non-Markovian drifting continuous-time random walk processes}

\author{Miquel Montero}
\email[E-mail: ]{miquel.montero@ub.edu}
\affiliation{Departament de F\'{\i}sica Fonamental, Universitat de
Barcelona, Diagonal 647, E-08028 Barcelona, Spain}
\author{Javier  Villarroel}
\email[E-mail: ]{javier@usal.es}
\affiliation{Facultad de Ciencias, Universidad de Salamanca, Plaza Merced s/n, E-37008 Salamanca, Spain}

\begin{abstract}
By appealing to  renewal theory we   determine   the equations
that the mean exit time of a continuous-time random walk with
drift  satisfies both when the present
 coincides with a jump instant or when it does not.  Particular attention
is paid to the corrections ensuing from the non-Markovian nature
 of the process.   We show that
 when drift and jumps have the same sign the relevant integral
 equations can be solved in closed form. The case when  holding
 times have the
 classical Erlang distribution is considered in detail.

\end{abstract}
\date{\today}
\maketitle

\section{Introduction}

 In this article we
 study exit times of  continuous-time random walks (CTRWs) with
 drift. By this we understand a
random process $X_t$ whose evolution in
 time
 can be thought of  as the result of the combined effect of a constant
drift and the occurrence of random
 jumps. Thus,  for $t\geqslant t_0$ we define
\begin{equation}
X_t=X_{t_0}+v(t-t_0)+S_{t}, \quad S_t=\sum_{n=1}^\infty J_n\
\theta(t-t_n), \label{process}
\end{equation}
where $\theta(u)=1$ for $u \geqslant 0$,
$t_n=\tau_1+\cdots+\tau_n$ are the jump  times,  $J_n$ the
jump magnitude and $v>0$ by convention. Note that as the notation
suggests, $t_0$ coincides with a jump.  The process $X_t$ may
represent physically the position of a diffusing particle, for
instance. All through this paper  we will assume that (i) the
sojourn times $\tau_n>0$
 are independent and identically distributed (i.i.d.) random variables with
probability density and cumulative distribution function (PDF and,
respectively, CDF) $\psi(t)$ and $\Psi(t)=  \int_0^t \psi(t')d
t'$; (ii) $\{J_n\}$ is a sequence of i.i.d. random variables  with
common PDF $h(\cdot)$; and  (iii), that  $ J_m$ is independent of $\tau_n $ for any $n,m$.

In absence of drift, and when  the holding times $\tau_n$ are exponentially distributed, $\tau_n\sim  \mathcal  E(\lambda)$ for some $\lambda >0$, the jump process  $S_t$ is  a  classical compound Poisson process (CPP)
characterized by having     independent increments  $S_{t+h}- S_t$; in addition the
    associated  ``arrival-process''   $N_t\equiv\sum_{n=1}^\infty \theta(t-t_n)$  is
  Poisson
 distributed:   $N_t \sim
\mathcal P(\lambda
 t)$. Therefore, drift-less CTRWs can be considered as a generalization of CPPs where the holding times of the  
processes $S_t$ are arbitrarily distributed. In  statistical physics  such  drift-less  CTRWs      have been widely
used  after the work of Montroll $\&$ Weiss \cite{mw65,w94} and
their relevant statistical properties, like correlation functions
and the behavior in the continuum limit,~\footnote{This limit corresponds to the assumption that the sojourn time between consecutive jumps and their characteristic size go to zero in an appropriate and coordinated way. We will show an explicit example of this limit in the Appendix.}
 a subject of intense interest
\cite{JT74,JT76,Go,Sch}.  Applications of CTRWs can be found in the study of transport in
disordered media (e.g., \cite{Sch,ms84,Weiss-porra,Margolin1}),
anomalous relaxation in polymer chains~\cite{Hu}, sandpile and earthquake 
 modeling (e.g., \cite{hs02,Me03}),   random networks \cite{Be},
 self-organized
criticality in granular systems~\cite{Bo},   scaling properties of  L\'evy walks~\cite{Ku}, electron tunneling~\cite{Gu}, transmission
 tomography~\cite{OW,OW2},
  distribution of matter in the universe~\cite{Os71} and changes  in  stock
markets due to unexpected catastrophes~\cite{m76}. More
recently,
 the use of CTRWs has been advocated to give a microscopic, tick-by-tick,
description of financial markets: see~\cite{mmw03,mmp05,mmpw06,mpmlmm05,s09}. A comprehensive review of CTRW applications in finance and economics is given in \cite{s06}.

Physically,
 the introduction of these general CTRWs stems from the fact that
in  many settings the exponential holding-time   assumption  may
be inadequate to describe the physical situation |see \cite{Me03,mmw03}.  Additional  motivation arises from
the observation that there is an extensive number of different physical systems that show some sort of anomalous diffusion |a subject of great interest from the viewpoint of statistical physics in the last years|, and that this anomalous behavior can well described by fractional Fokker-Planck equations obtained after imposing the continuum limit on CTRW models~\cite{
Metzler,s07,m07}.

The further addition of the drift term to a CTRW, as we do in Eq.~(\ref{process}), is a natural and significant incorporation. The resulting process |which can be viewed as the discrete  analogue of a (fractional) diffusion with drift| is known for playing  a fundamental   role  in 
the modeling of
the cash flow at an insurance company \cite{j03,zyl10} and, more recently,  it has been shown that it 
also rules the rate of   energy dissipation in  nonlinear optical
fibers~\cite{VM1}. Indeed, present and forthcoming results  can be also of interest in
   transport in amorphous media~\cite{Ba},
models of  decision
  and response time in psychology~\cite{Sm} and neuron dynamics~\cite{Sm2}. 

In all these scenarios one is faced with the basic problem of determining the first-passage time for a CTRW with drift, a question that has been theoretically considered in the past~\cite{m87,m92,c97,cmc97,rd00a,rd00b,MB,is07}. The usual approach taken there entails the computation of transition and first-passage time distributions of the process, e.g. \cite{m92,MB}, and the results are typically obtained under the assumption of the continuum limit~\cite{rd00a,rd00b,MB,is07}. 
This procedure has incontestable advantages for obtaining the leading-order behavior, but it is not adequate for a detailed analysis of the statistical properties of the process at the inter-jump timescale, as in~\cite{mmw03,mpmlmm05,
j03,zyl10,VM1}.

Motivated by the above,
here we pose the problem of evaluating  the mean exit time from
the  interval 
$(0,b)$ of a drifting CTRW $X_t$ when {\it the only available information is the present state\/} $X_r=x$, where $r \geqslant 0$ is the present chronological time 
|note that, by
adjusting the time clock and spatial scale, the results carry
over to any interval $(a,b)$ and initial time $t_0$. Here
$X_t$ is given by Eq.~(\ref{process}) where 
$S_t$ is the jump part and $N_t$
the associated counting process (a renewal
    process).    In the drift-less case previous work in this regard includes that
of \cite{mmp05,mmpw06} where  a linear integral equation for the
mean escape time after a jump off a given interval is derived. We
note however that these
 results  do not cover a generic situation. Indeed, while  for
CPPs (as  for  the  general  L\'evy processes) Markov property
implies that results derived starting at a jump time carry over
to arbitrary present, no such inference is possible for a generic
CTRW  due to its non-Markovian nature. Thus, escape times \it
depend on the actual state and available information \rm and
hence  the question as to how to generalize the former results to
\it general present time \rm $r$ appears naturally. We remark that
 implications  ensuing from the lack of
 Markovianess
 have been      ignored at large in the literature, a gap that we
  have intended  to fill in |see \cite{VM2}.   In particular, it  remains
 an open
question  to what extent  dropping   the assumption  that ``the
present is a jump
   time'' affects  the relevant probabilities.
 Here  we address these issues      and
generalize the  results of  \cite{mmp05,mmpw06} in a twofold way
by assuming  that   (i) a drift $vt>0$ operates on the system and
(ii) the present $r$ is an arbitrary   time,    not necessarily a
jump instant, and the observer \it  has  knowledge of the present,
but not of  the history, of the system\rm.

The interest of this problem goes far beyond the purely academic
 since  such a situation  may  appear in
several different physical contexts. For example, one might be
interested in predicting the mean time  for an
insurance/financial  company to go bankrupt   from the knowledge
of just the actual  company budget, i.e.,  when the information
regarding the company's past performance has not been disclosed.
A second example is provided by the  study of   the distribution
of inhomogeneities in an optical fiber;   it was  found  \cite{VM1}    that the signal's energy  amplitude at a point $t$
($t$ is the spatial variable in this setting)  involves  a CTRW
with drift $X_t$; in this context~\footnote{This situation where
$t$ represents the {\it the space variable } may occur in
different physical systems,
 with  $X$ standing  for some physical observable of interest,
like the energy. In  this connection, $r$ might represent the
location of a detector, or a sensitive part of the appliance at
which measures are taken.}  one typically knows only the value of
the  energy at the observation point $r$, not on the whole fiber.
More generally the approach will be relevant in situations where
  either the elapsed time between events  is ``large''  (it might be as large as
  years, in a context of catastrophes observation)
   or when the event's times $t_n$ are not physically
measurable  observables and only mean escape times are. (Note that
both the mean exit time and initial time are typically   {\it
macroscopic\/} magnitudes.)

 The article is structured as follows. In Section~\ref{S2} we show how the solution to these problems involves ideas drawn from renewal theory
 and solve the simpler case when $v=0$. The case $v>0$ is
 considered in the next sections where it is found that 
 key properties of the obtained equation depend on the sign of the
 jumps. In Section~\ref{S_up} it is shown  that if this sign is positive the solution can be
 given in closed form by Laplace transformation |cf.
 Eqs.~(\ref{FF_sol}) and~(\ref{hattime}). For the case  when drift and jumps have opposite signs
  we find integral equations  that  the relevant objects satisfy, see Section~\ref{S_down}, but no closed
 solution can be given in a general  situation. Section~\ref{S_two-sided}
 addresses the most general scenario in which
 jumps  $J_n$   can   take both signs.   Solvable cases are
 discussed there.

  In all cases we exemplify our results by considering the particular instance when
    sojourn times have
Erlang distribution, $\mathcal Er(\lambda,2)$. 
 $\mathcal Er(\lambda,n)$   corresponds to having a sum of $n$ independent exponential
variables and hence generalizes the exponential density in a
natural way, 
\begin{equation}
\psi(t)  = \frac{\lambda^n t^{n-1}}{(n-1)!} e^{-\lambda t},\, n\in \Bbb N,\, t\geqslant 0,  
\label{psi_E2}
\end{equation}
while it maintains an adequate capability to fit
measured data. From a physical perspective   these facts make
 this density   a natural   candidate to describe  multi-component systems  which
operate  only when several independent, exponentially distributed
operations have been completed or whenever  there is  a hidden
Poissonian flux of information  and jumps only appear  as the
outcome of  two or more  consecutive arrivals.
 This explains the interest that
 it has drawn in the field  of information traffic
\cite{sb00,fc02}.   Similarly the appearance of this distribution
to model transaction orders in financial markets can also be
expected since
 it takes, at least, two arrivals (buy and sell orders) for  a transaction
to be completed.   For further applications to ruin problems and
insurance see  \cite{dh01,lg04}.

\section{\label{S2}The problem}

Recall that we  aim to   study exit times of a  drifting CTRW
$X_t$
 given the present state $ X_r=x$.  To this  end    let $r+\mathbf{t}_b^{x,r}$ be the first time past
$r$ at which $X_t$ exits $(0,b)$;   $T_{b}(x,r)$ be its expected
value: $T_{b}(x,r) =\EE[\mathbf{t}_b^{x,r}]$; finally let
$\TT_{b}(x) $ denote   the mean exit time   off $(0,b)$     after
a jump $t_n$ (loosely one has $\TT_{b}(x)=T_{b}(x,t_n), n=0,\dots,
\infty$; however the relation between both   quantities is not
trivial, as we see below in Section~\ref{S_up}). Note also that here and
elsewhere we use $\Bbb E[\cdot]$ to denote expectation.

   In the exponential Markov case    $
\mathbf{t}_b^{x,r}$ is independent of $r$, $T_{b}(x,r)=\TT_{b}(x)$
and it \it only \rm remains to formulate (and solve) the equation
that  this object satisfies. However this situation no longer
holds in the generic, non-Markovian case where    $
\mathbf{t}_b^{x,r}$ does depend on $r$. We find (see below and
Section~\ref{S_up}) that  the relation between $T_{b}(x,r)$ and
$\TT_{b}(x)$ involves     the  ``excess life'' $E_r\equiv
t_{N_r+1}-r$, or time elapsed   until  the next arrival    occurs.
We  now  sketch classical renewal theory (see \cite{kt81,c65}
 and  \cite{Go}) that
   shows how to construct  the CDF    $\Phi(t|r)\equiv
\Bbb P\{E_r\leqslant t\}$ of $E_r$.

   Let  $m(t)\equiv \EE[N_t]$   be  the  mean
number of jumps up to $t$:  the renewal function. It satisfies
the integral renewal equation
\begin{equation}
m(t)=  \Psi(t) +\int_0^t  m(t-t')\psi(t')d  t'. \label{re}%
\end{equation}
 Then, by using the total
probability theorem   it can be proved  that
\begin{equation}
\Phi(t|r) = \int^{r+ t }_{r}\left[1- \Psi(r+ t-t')\right]d  m(t').
\label{repp}%
\end{equation}  Upon solution of the  above
  integral  equations we obtain $\Phi(t|r)$. Actually, they  can be solved
  with all generality
  by recourse to   Laplace
  transformation.   Let $\hat g(s)$ be   the Laplace
transform of a function $g(t)$ so that
\begin{equation}  g(t) =\frac{1}{2\pi i}
\int_{c-i\infty}^{c+i\infty} e^{ st} \hat g(s)  d  s,\ c>0.
\label{def_m}%
\end{equation}
  Then,    Eqs.  (\ref{re}) and (\ref{repp}) allow  to  recover the
  distribution of $E_r$ in closed form via
\begin{subequations}
\label{hat_phi} 
\begin{eqnarray}
\hat m(s)&=&\frac{1}{s} \frac{\hat\psi(s)}{1-\hat \psi(s)}, \\   
\hat \phi(s|r)&=&  e^{s  r } \left[1-\hat \psi(s)\right] \int_r^\infty e^{-sl} \hat m(l)l d  l, 
\end{eqnarray}
\end{subequations}
where   $ \hat m(s) $ and $  \hat \phi(s|r)$ are the Laplace transforms of  $m(t)$  and $\phi(t|r) \equiv \p_t \Phi(t|r)$.

 If $v=0$ these expressions can be used to relate $T_{b}(x,r)  $ and
 $\TT_{b}(x)$. Indeed  let $t_{N_r}\eq t_n$, say,  be the ``last'' jump time
   and  $E_r^-  \equiv r- t_n$     the time elapsed
  from  $t_n$     to the present. Then, with $E_r^+\equiv
 E_r$ one  obviously has that
 $E_r^-+E_r^+=t_{n+1}-t_n\equiv\tau_{n+1}$, and that the exit time right
after $t_n$ is that after $r$, $\mathbf{t}_b^{x,r}$,  plus $E_r^-$, and hence
 \begin{equation}\TT_{b}(x)=\EE [E_r^-]+
\EE[\mathbf{t}_b^{x,r}]=\EE[\tau_{n+1}]-\EE[E_r^+]+
T_{b}(x,r).\label{v=0}%
\end{equation} Thus,
 $T_{b}(x,r)
 $   follows    adding a   correction     term to $\TT_{b}(x)$ which  depends only on $r$. Finally $\TT_{b}(x)$
 is obtained by solving  a linear
 integral equation |see \cite{mmp05,mmpw06}.
Unfortunately when $v\ne 0$ this simple argument  fails  as
 then \it  knowledge of the present position does not entail its   knowledge
  at $t_n$\rm.
 In the next sections   we   derive  the relevant correction
 to
the mean
  exit time. This correction depends now in all parameters $r$, $x$ and $b$ |see Eqs.~(\ref{Tb_up}) and~(\ref{Tb_down})  below.

\section{\label{S_up}Jump process with favorable drift}
In this section we consider the case when both drift and jumps
have a positive sign, i.e. when $X_t$ is increasing. As a result, the process can only leave the interval through the upper
boundary $b$. Let us assume that at time $t=r$ the system is in $X_r=x\in(0,b)$, and that the excess life $E_r$ is known in advance, $E_r=l$. If the excess life is longer than $\varrho \equiv {b-x\over v}< l$, the drift will drive the process
out of the region at time $r+\varrho$, before the next jump takes place at $t_{N_r+1}$. Conversely, if $l \leqslant \varrho$ at least a jump of size $J_{N_r+1}=u$ will occur prior to exiting the interval. Note that just before $t_{N_r+1}$ the process $X_t$ is no longer at $x$, but at $x+v l$. Therefore, two possible scenarios appear: either
the jump size is larger than the remaining distance up to the upper boundary, $u\geqslant b-x-v l$, and the process leaves the
interval at $r+E_r$, 
or it does not. In the latter case the problem renews from $t=r+l$, $X_{r+l}=x+v
l+u$, so the mean escape time will be increased by an amount $\TT_b\left(x+v
l+u\right)$. It can be proven that  these considerations imply
that $T_b(x,r)$ must satisfy |recall that $\psi(\cdot)$, $h(\cdot)$ are the
waiting-time and, respectively, jump PDFs| the following equation
\begin{eqnarray}
T_b(x,r)=\left[1-\Phi\left(\varrho|r\right)\right] \varrho +
\int^{\varrho}_{0} l
\phi(l|r) \int^{\infty}_{b-x-v l} h(u) d  ud  l\nonumber \\
+\int^{\varrho}_{0} \phi(l|r) \int^{b-x-v l}_{0} h(u)
\left[l+\TT_b\left(x+vl+u\right)\right]d  ud  l,
\label{Tb_up0}
\end{eqnarray} an expression that relates $T_b(x,r)$ and $\TT_b(x)$.
Note that since $\Phi\left(l|r\right)$   depends on   $r$ so it
does the mean time $T_b(x,r)$. Finally, after some rearrangement, Eq.~(\ref{Tb_up0}) can be conveniently written as
\begin{eqnarray}
T_b(x,r)= \int^{\varrho}_{0}\left[1-\Phi\left(l|r\right)\right]d
l \nonumber \\+ \frac{1}{v}\int^{v \varrho}_{0}
 \phi\left(\left.\varrho-\frac{z}{v}\right|r\right) \int^{z}_{0} h(u)
\TT_b\left(b-z+u\right)d  u  d  z. \label{Tb_up}%
\end{eqnarray}

 Similarly, by
 letting 
 $r\to 0$, we find that
 $\TT_b(x)$ must satisfy
\begin{eqnarray}
\TT_b(x)= \int^{\varrho}_{0}\left[1-\Psi(l)\right]d  l  \nonumber \\
+\frac{1}{v}\int^{v \varrho}_{0}
\psi\left(\varrho-\frac{z}{v}\right) \int^{z}_{0} h(u)
\TT_b\left(b-z+u\right)d  ud  z, \label{TTb_up}%
\end{eqnarray}
$x\in(0,b)$. Equation~(\ref{Tb_up})  along with~(\ref{TTb_up})  allow us to solve
the posed problem. The second of these
  defines an  integral equation for $\TT_b(x)$ which, if $v=0$,
 reduces upon appropriate change to that  of  \cite{mmp05,mmpw06}.
  It is remarkable that Eq.~(\ref{TTb_up}) \it can be solved in a fully
explicit way. \rm   To this end we define the    allied object
$\FF(y)$, 
as the solution of the following integral
equation
\begin{eqnarray}
\FF(y)= \int^{y/v}_{0}\left[1-\Psi(l)\right]d  l  \nonumber \\
+\frac{1}{v}\int^{y}_{0} \psi\left(\frac{y-z}{v}\right)
\int^{z}_{0} h(u) \FF\left(z-u\right)d  ud  z, \label{FF_eq}%
\end{eqnarray}
for $y\in \RR^{+}$. Then it follows that $\TT_b(x) =\FF(b-x)$, for
$x\in (0,b)$.  We note further that   taking a   Laplace
transformation in   Eq.~(\ref{FF_eq})  we    find that
\begin{equation}
\hat{\FF}(s)= \frac{1}{v s^2}\frac{1-\hat{\psi}(s
v)}{1-\hat{\psi}(s v)\hat{h}(s)}, \label{FF_sol}
\end{equation}
where  $\hat{\FF}(s)\equiv \int^{\infty}_{0} \FF(y) e^{-s y}d y$,
$\hat{h}(s)\equiv \int^{\infty}_{0} h(u) e^{-s u}d  u$, and
$\hat{\psi}(s)\equiv \int^{\infty}_{0} \psi(s) e^{-s t}d  t$. Here
$s=s_R+is_I$ is complex and $s_R\geqslant 0$. Further, we also have
\begin{equation}
\hat{J}(s|r) =  \hat{\FF}(s)-   \frac{ 1- \hat{h}(s) }{v
s^2\left[1-\hat{\psi}(s v)\hat{h}(s)\right]}\left[\hat{\phi}(s
v|r)-\hat{\psi}(s v )\right],\label{hattime}
\end{equation}
where the correction to the mean time after a jump is clearly
displayed and again, for convenience, we defined $J(y|r)\eq T_b(b-y|r)$, and extended (\ref{Tb_up}) to $y\geqslant 0$. Thus ${\FF}(y) $ and $J(y|r)$  can be recovered by
Laplace inversion, cf. Eq.~(\ref{def_m}), and $\TT_b(x)$ and $T_b(x|r)$ will eventually follow.

Equation~(\ref{hattime}) has several mathematical limits of applied interest. We first consider the situation when $v$ is small. Let $\mu\equiv~\EE[\tau_{n+1}]$, $\mu_r\equiv~\EE[E_{r}]$; then, using   that $  \hat{\psi}(s v )=1-\mu sv
+O(v^2)$ and so forth we see that for  small $v$, $\hat{ J}(s|r)$
has an expansion in powers of $v$ as
 \begin{equation}
\hat{J}(s|r) =  \hat{\FF}(s) - \frac{\mu-\mu_r}{s}  -    \frac{
\hat{h}(s) }{ 1- \hat{h}(s) }   (\mu-\mu_r)\mu  v + O(v^2),
\end{equation}
which implies, in particular, Eq.~(\ref{v=0}). Another interesting
case is obtained letting $r\to\infty$: the steady-state solution.
This limit is relevant since it  can be associated to a
situation  in which \it the only information available to the
observer is the present value of the stochastic process, not even
the starting point\rm. Recalling that by the renewal theorem
$\un{t\to\infty}\lim m(t)/t = 1/\mu$,  Eqs.~(\ref{hat_phi}) and~(\ref{hattime}) yield
 \begin{eqnarray} \hat\phi (s|\infty)
&=&  \frac{1-\hat \psi(s)}{s\mu},  \\
\hat{J}(s|\infty) &=&
\hat{\FF}(s)- \frac{ 1- \hat{h}(s) }{v^2 s^3 \mu}
\frac{1-(1+s v\mu)\hat \psi(sv)}{ 1-\hat{\psi}(s
v)\hat{h}(s)}.
\end{eqnarray}

 We  illustrate the ideas above by detailing the case when jump magnitudes
  have
an exponential distribution
  $h(u)=\gamma e^{-\gamma u}$, where $\gamma>0$ is a real parameter, and  sojourn times
 an Erlang distribution $\tau_n\sim\mathcal Er(\lambda,2)$.
Hence
\begin{equation}
\hat{h}(s)= \frac{\gamma}{\gamma +s},\quad  \hat{\psi}(s)  = \frac{\lambda^2}{(\lambda +s)^2}. 
\label{h_hat}
\end{equation}
In this case Eq.~(\ref{FF_eq})  yields that $\hat{\FF}(s)$ is the
following rational function
 \begin{equation}
\hat{\FF}(s)=\frac{1}{v s^2}\frac{2\lambda \gamma v+ v(2 \lambda
+\gamma v) s + v^2 s^2}{\lambda (\lambda + 2\gamma v) + v(2
\lambda +\gamma v) s+v^2 s^2}.
 \end{equation}
 Hence, the mean  escape  time of $(0,b)$ after a jump is
  \begin{eqnarray}
  \TT_b(x)= \frac{2 \gamma v}{\lambda + 2 \gamma v} \varrho \nonumber\\
+\frac{\lambda^2}{\zrp-\zrm}\left[\frac{1-e^{-\zrm \varrho}}{\zrm^2}-\frac{1-e^{-\zrp
\varrho}}{\zrp^2}\right],
\end{eqnarray}
where we recall that $v \varrho=b-x$ is   the  initial
      distance to the boundary $b$
  and
\begin{equation}
z_{\pm}\equiv
 \lambda +\frac{\gamma v}{2}\left(1\pm\sqrt{1-\frac{4\lambda}{\gamma
v}}\right),
\end{equation}
with $\mbox{Re}[z_{\pm}]>0$.

\begin{figure}[htb]
{
\includegraphics[width=0.45\textwidth,keepaspectratio=true]{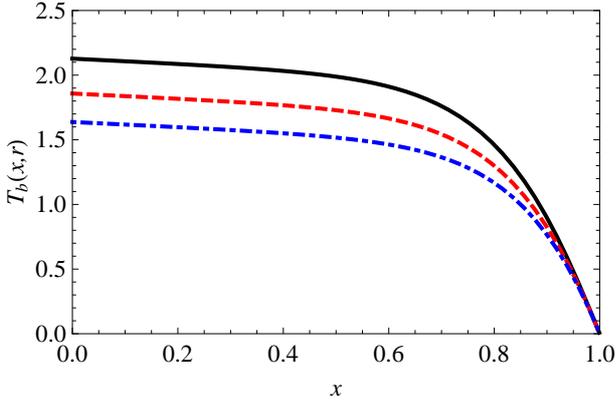}
} \caption{(Color online)  Sample representations of $T_b(x,r)$ for different values of $r$: $r=0$, i.e. $\TT_b(x)$, solid (black) line, $r=0.4$,  dashed (red) line and $r=10.0$, dot-dashed (blue) line. The rest of the parameters were chosen as follows: $b=1.0$, $v=0.1$, $\lambda=1.0$ and $\gamma=0.1$.} \label{E2_up}
\end{figure}

The evaluation of  the mean exit time  starting at $r$,
$T_b(x,r)$, involves the renewal function and excess-life
distribution. We first obtain    from     Eq.~(\ref{hat_phi})
\begin{equation}
\hat m(s)=\frac{\lambda^2}{s^2(2\lambda +s)}, \quad  \hat
\phi(s|r)=\hat \psi(s)+ \frac{  (1-e^{-2 \lambda r})s\lambda }{ 2
(\lambda +s)^2}.
\end{equation}
By inversion we get $m(t)=  \big(2 \lambda t + e^{-2 \lambda
t}-1\big)/4$ and
\begin{equation}  \ \Phi(t|r)  
=
1-e^{-\lambda t} \left[ 1+  \left(\frac{ 1+{e^{-2 \lambda r
}}}{2}\right)\lambda t \right].
\end{equation}
Then by using  Eq.~(\ref{hattime})  we find at last that
\begin{eqnarray}  T_b(x,r) =\TT_b(x) \nonumber 
\\
-  \frac{\lambda}{2}\frac{1-e^{-2\lambda r}}{\zrp-\zrm}\left[\frac{1-e^{-\zrm \varrho}}{\zrm}-\frac{1-e^{-\zrp \varrho}}{\zrp}\right].
\end{eqnarray}
Plots of this function in terms of $x$ are given in Fig.~\ref{E2_up} for several values of $r$ and a certain choice of the rest of parameters.

Finally note how, in particular,  if $v=0$     then  
\begin{eqnarray}
\TT_b(x)=\frac{2}{\lambda}\left[1+\gamma (b-x)\right] &=&  \EE[\tau_n]
\left(1+   \frac{b-x}{\EE[J_n]}\right), \\
T_{b}(x,r)&=&\TT_b(x)
 -\frac{1- e^{-2\lambda  r}}{2\la}.\end{eqnarray}

\section{\label{S_down}The case of opposite drift and jumps}

 We now consider the  case when  the sign  of the drift $vt>0$  is
opposite to that of the   jumps.
  In this case the process can leave the
interval through both the upper and the lower boundaries: the
drift pushes steadily the process up, whereas the jumps threaten
the system with a downside exit.   The   resulting process  is
 a prototype  model in  risk management to describe   the
dynamics of the cahsflow  $X_t$ at an insurance company  under
the assumption  that premiums are received at a constant rate
$v>0$ and that  the company incurs in losses $J_n<0$
  from
claims reported at times $t_n, n=1,\dots\infty$   (the
Cramer-Lundberg model).

 As before, we analyze the evolution starting  at $t=r$ with $X_r=x$. Then if
$E_r=l$, $l > \varrho$,
  the drift will drive the
process out of the region through the {\it upper boundary\/} at time
$r+\varrho$. Otherwise  ($l \leqslant \varrho$) at least a jump
$J_{N_r+1}=-u$,  say, occurs at time $t_{N_r+1}$  before escape, and two possible scenarios appear depending on the relative magnitudes of the jump size and the location of the process right before the jump, $X_t=x+v l$: for $u\geqslant x+v l$ the
process will leave the interval through the {\it lower boundary\/} at $r+E_r$, but when $u<x+v l$ the process after the jump will remain inside the interval, $X_{r+l}=x+v l-u$, the exit problem will start afresh, and the mean escape time will be increased by $\TT_b\left(x+v l-u\right)$.  Again these considerations imply
that $T_b(x,r)$ and $\TT_b(x)$ must satisfy for $0< x< b$
\begin{eqnarray}
T_b(x,r)=\left[1-\Phi\left(\varrho|r\right)\right] \varrho +
\int^{\varrho}_{0} l
\phi(l|r) \int^{\infty}_{x+v l} h(u) d  ud  l \nonumber \\
+\int^{\varrho}_{0} \phi(l|r) \int^{x+vl}_{0} h(u)
 \left[l+\TT_b\left(x+vl-u\right)\right]d  ud  l\nonumber \\
=
\int^{\varrho}_{0}\left[1-\Phi\left(l|r\right)\right]d  l \nonumber \\
+
\frac{1}{v}\int^{b}_{x}
\phi\left(\left.\frac{z-x}{v}\right|r\right) \int^{z}_{0} h(u)
\TT_b\left(z-u\right)d  ud  z,\nonumber \\ \label{Tb_down}  
\end{eqnarray}
\begin{eqnarray}
\TT_b(x)= \int^{\varrho}_{0}\left[1-\Psi(l)\right]d  l
\nonumber \\
+\frac{1}{v} \int^{b}_{x} \psi\left(\frac{z-x}{v}\right)
\int^{z}_{0} h(u) \TT_b\left(z-u\right)d  ud  z. \nonumber \\ \label{TTb_down}
\end{eqnarray}

 Hence
 $T_b(x,r)$ follows in terms of quadratures also in this case,  given $\TT_b(x)$.
  Unfortunately, unlike  what happens   for  the case    considered in the previous Section,
  Eq.~(\ref{TTb_down})   can not be  solved in  closed form
   for   arbitrary PDFs $\psi(\cdot)$ and
$h(\cdot)$.    Further progress can be made for Erlang times, $\mathcal Er(\lambda,2)$. 
Indeed, in this case (\ref{TTb_down}) reads
\begin{eqnarray}
\TT_b(x)= \frac{2-\left(2+\lambda \varrho\right)e^{-\lambda
\varrho}}{\lambda} \nonumber \\+\frac{\lambda^2}{v^2} \int^{b}_{x}
\left(z-x\right)e^{-\lambda \left(z-x\right)/v} \int^{z}_{0} h(u)
\TT_b\left(z-u\right)d  ud  z,
\end{eqnarray}
 and hence, by differentiation we find that
$ \TT_b(x)$, for $0< x< b$,  also satisfies the following
integral-differential equation:
\begin{eqnarray}
\TT_b''(x)-\frac{2\lambda}{v}\TT_b'(x)+\frac{\lambda^2}{v^2}\TT_b(x)=
\frac{2\lambda}{v^2} \nonumber \\+\frac{\lambda^2}{v^2}\int^{x}_{0}h(u) \TT_b\left(x-u\right)d  u, 
\label{IDb_down}
\end{eqnarray}
subject to the following boundary conditions: 
\begin{equation}
\lim_{x\to b}\TT_b(x)=0,\mbox{ and }\lim_{x\to b}\TT_b'(x)=-1/v.
\label{boundary}
\end{equation}

We first consider a general solution to this equation extending
it   to the full  real axis, so we will drop the subscript in $\TT_b(x)$.
    We find a solution by Laplace transformation as
\begin{eqnarray}
\hat{ \TT }(s)= \frac{2\lambda/s +   (B v-2 \lambda A)v  + A
v^2s}{\lambda^2[1-\hat{h}(s)]-2\lambda v s+ v^2s^2},
\end{eqnarray}
where   $A $ and  $B$ are 
$\TT(0)$ and $\TT'(0)$ respectively. By
inversion, cf. Eq.~(\ref{def_m}), $\TT(x)$ follows in terms of
$A$ and $B$. By requiring~(\ref{boundary}) 
we obtain a linear algebraic system for $A$ and $B$,
which upon solution yields $\TT_b(x)$ in closed form.

 To be specific we consider the case when jumps are also exponentially
distributed: $h(x)=\gamma e^{-\gamma x}$. Then we have~(\ref{h_hat}) 
and  $\hat{\TT}(s)$ is the rational
function
\begin{eqnarray}
\hat{\TT}(s)= \frac{(\gamma+s)[2\lambda + (B v-2 \lambda A)v s+ A
v^2s^2]}{s^2[\lambda(\lambda-2 \gamma v)+(\gamma v-2\lambda) v s+
v^2s^2]}. \label{LUps}
\end{eqnarray}
Upon re-scale  of constants the inverse Laplace transform
of~(\ref{LUps}) reads

\begin{eqnarray}
\left(\frac{\lambda}{2}-\gamma v\right)\TT(x)\nonumber \\
=\tilde A  +  \gamma   x +
 \left[\frac{\lambda}{v^2}(\lambda-\lambda
\tilde A+   v\tilde  B) -\frac{\tilde B}{2}
\xi_{-}\right]\frac{e^{\xi_+ x}-1}{\xi_+(\xi_+-\xi_-)} \nonumber \\
- \left[\frac{\lambda}{v^2}(\lambda-\lambda\tilde  A+   v
\tilde B) -\frac{\tilde B}{2}\xi_{+}\right]\frac{e^{\xi_-
x}-1}{\xi_-(\xi_+-\xi_-)},\nonumber \\
\end{eqnarray}
with
\begin{equation}
\xi_{\pm}\equiv
\frac{\lambda}{v}-\frac{\gamma}{2}\pm\frac{\gamma}{2}\sqrt{1+\frac{4\lambda}{\gamma
v}}.\end{equation}

\begin{figure}[htbp]
{
\begin{tabular}{rl}
(a)&\includegraphics[width=0.45\textwidth,keepaspectratio=true]{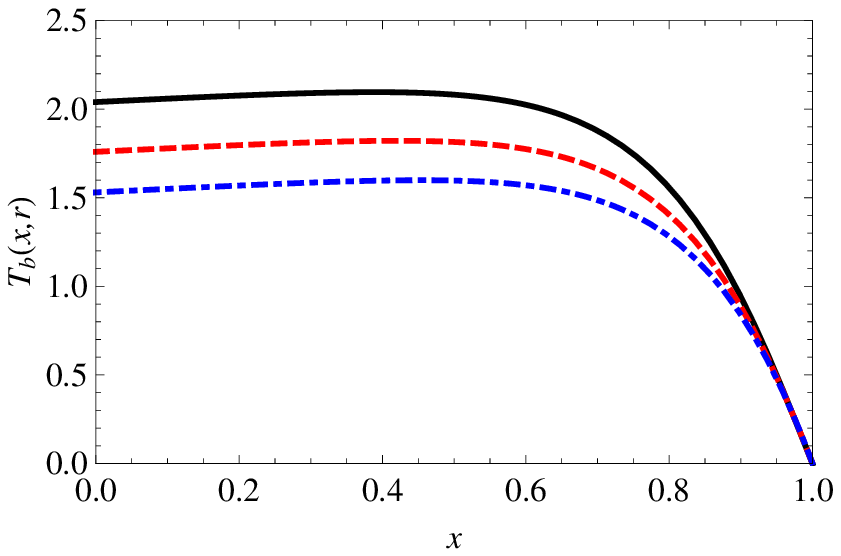}\\
(b)&\includegraphics[width=0.45\textwidth,keepaspectratio=true]{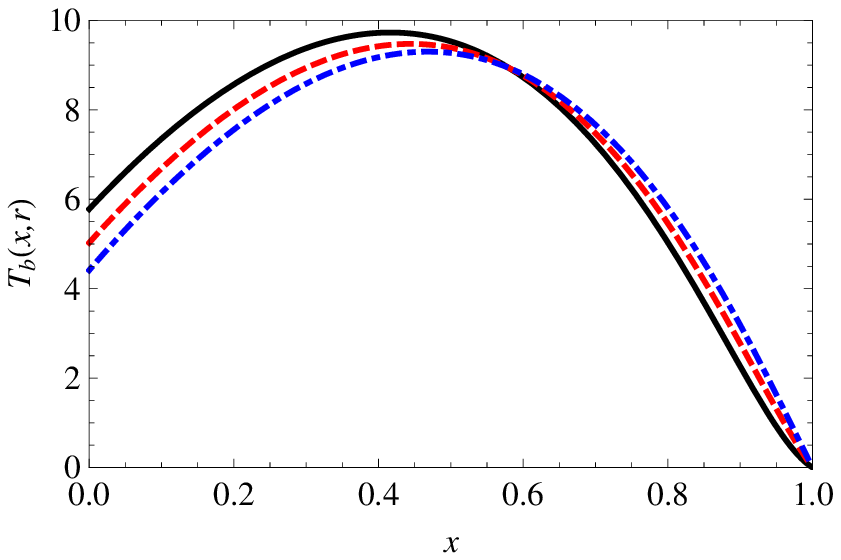}\\
\end{tabular}
}
\caption{(Color online) Sample representations of $T_b(x,r)$ for different values of $r$: $r=0$, i.e. $\TT_b(x)$, solid (black) line, $r=0.4$,  dashed (red) line and $r=10.0$, dot-dashed (blue) line. The rest of the parameters were chosen as follows, $b=1.0$, $v=0.1$, $\lambda=1.0$ and (a) $\gamma=0.1$; (b) $\gamma=4.0$. 
We see how, when drift and jumps have a different sign, an interesting behavior may appear in $T_b(x,r)$.
In the upper panel we observe that this function is no longer decreasing with $x$ and has a maximum in the interior of the interval, cf. Fig.~\ref{E2_up}. In the lower panel we observe how a cross-over phenomenon may eventually take place reflecting the fact that $\TT_b(x)$ need not be greater than $T_b(x, r)$: this behavior can be traced to the fact that as $r$ increases the jump probability increases, which results in a reduction of the escape time if the process is in the vicinity of the lower boundary but in an increment when close to the upper boundary. The maxima position moves toward greater values of $x$ for larger choices of $r$, in both cases.
} \label{E2_down}
\end{figure}

Unfortunately    the final  expressions for $\tilde A$ and $\tilde B$  after imposing~(\ref{boundary}) are not
  very illuminating so
 we   do not  transcribe them here. Sampling values for different
parameter specifications can be found in Fig.~\ref{E2_down}.

The limit  $b\to \infty$ is interesting as $\TT_\infty(x)$ gives
the probability that $X_t$ ever hits $0$. This corresponds to the
classical ruin probability in an insurance context.   It turns
out that $\TT_\infty(x)$ can be determined in a direct way that
avoids  solving the aforementioned linear system. Without
proof~\footnote{We elaborate on a similar expression in the next
section.} we note that if $\lambda
> 2 \gamma
 v $ then
\begin{equation}
\TT_\infty(x)=\frac{2(1+\gamma x)}{\lambda-2 \gamma v}= \left(  1+\gamma
x  \right)\Big/\left( \frac{1}{ \Bbb E[\tau_n]}- \frac{v}{  \Bbb E[J_n]}
\right),
\end{equation}
while $\TT_\infty(x)=\infty$ otherwise. Once $\TT_b(x)$ is
known, $T_b(x,t)$ follows again by integration |see Fig.~\ref{E2_down}.

\section{\label{S_two-sided}Two-sided jump process with drift}
We finally consider the general case corresponding to a jump
process where $J_n$ can take both signs and hence can    exit
$(0,b)$ through either of the boundaries.   The relevant analysis
is similar to that of the last section if one incorporates  the
possibility of an upper exit due to a jump: If $E_r=l>\varrho$ the drift drives directly the process through the upper boundary. If $l\leqslant \varrho$ at least a jump occurs before the exit of the process. When the size $u$ of this jump is positive and $u\geqslant b-x-vl$ the process leaves the interval at time $r+E_r$ through the upper boundary, when $u\leqslant -(x+vl)<0$ the exit takes place through the lower boundary, otherwise the process remains inside the interval $(0,b)$ and the problem restarts.
Skipping minor details
we obtain that $\TT_b(x)$ and $T_b(x,r)$ satisfy respectively
\begin{eqnarray}
T_b(x,r)=
\int^{\varrho}_{0}\left[1-\Phi\left(l|r\right)\right]d  l \nonumber \\+
\frac{1}{v}\int^{b}_{x}
\phi\left(\left.\frac{z-x}{v}\right|r\right) \int^{b-z}_{-z} h(u)
\TT_b\left(z+u\right)d  ud  z, 
\end{eqnarray}
and
\begin{eqnarray}
\TT_b(x)= \int^{\varrho}_{0}\left[1-\Psi\left(l\right)\right]d
l \nonumber \\+ \frac{1}{v}\int^{b}_{x} \psi\left(\frac{z-x}{v}\right)
\int^{b-z}_{-z} h(u) \TT_b\left(z+u\right)dud  z.  
\end{eqnarray}
   In  a general situation, the latter integral  equation is not
   solvable  in  closed form. To gain some insight  we use the
    decomposition
$ h(x)= q h_+(x)  + p  h_-(x)$ where $p$ is the probability  that
a given jump be  negative, $q\equiv 1-p$  and
 $h_{\pm}(\cdot)$  are the jump PDF
in the positive/negative regions, i.e.  $h_+(x)\equiv h(x)
\theta(x) \big/ q$, say.

Note that if $h_-(x)=0, \fo x\leqslant 0$, we recover the case considered
in Section~\ref{S_up}, solvable  via Laplace
    transformation. It turns out that we can still \it construct an
analytic closed solution \rm in
    the more general case when $h_-(x)  $  vanishes \it only \rm on $(-b,0)$ |but not on
    $(-\infty,-b]$.
     In such
a case a negative jump will  drive  the process out of the
interval $(0,b)$ through the lower boundary. Thus  $p$ is related
to    the ruin risk in a economic scenario or to  the breakdown
probability in a physical system.
 The  equation  for
 $\TT_b(x)$ reads in this case
\begin{eqnarray}
\TT_b(x)= \int^{\varrho}_{0}\left[1-\Psi\left(l\right)\right]d
l \nonumber\\+ \frac{q}{v}\int^{b}_{x} \psi\left(\frac{z-x}{v}\right)
\int^{b-z}_{0} h_+(u) \TT_b\left(z+u\right)d  ud  z,
\label{TTb_ruin}
\end{eqnarray}
while a similar expression, with $\phi(\cdot|r)$ in place of $\psi(\cdot)$, holds for $T_b(x,r)$.
Note that these equations are independent of the form of $h_-(\cdot)$
and  apart from the factor $q$ in front of the second term they
resemble Eqs.~(\ref{Tb_up}) and~(\ref{TTb_up}); it follows that
we can resort to the same technique used in Section~\ref{S_up}:
We consider
\begin{eqnarray}
\FF(y)= \int^{y/v}_{0}\left[1-\Psi(l)\right]d  l  \nonumber \\+
\frac{q}{v}\int^{y}_{0} \psi\left(\frac{y-z}{v}\right)
\int^{z}_{0} h(u) \FF\left(z-u\right)d  ud  z, 
\end{eqnarray}
for $y\in \RR^{+}$ which is again solvable  by means of a Laplace
transform; then $\TT_b(x)$ follows from  $\TT_b(x) =\FF(b-x)$  for
$x\in (0,b)$. With the previous selection for $h(\cdot)$ and
$\psi(\cdot)$, Eq.~(\ref{h_hat}), and in terms of    $z=sv $ 
we find

    \begin{eqnarray}
\hat{\FF}(s)=\frac{v}{z} \frac{2\lambda \gamma v  +  (2 \lambda+\gamma v ) z +  z^2}{p \lambda^2 \gamma v   + \lambda  (\lambda +2\gamma v  )  z+   (2 \lambda +\gamma v  )  z^2+z^3  }.\label{FF}
\end{eqnarray}
We first consider the case when   $\lambda
  =\gamma v$. Under this assumption  $\hat{\FF}(s)$ has poles at
$z=z_j\equiv \lambda \left(q^{1/3} e^{2\pi ij/3}-1\right)$, $j\in\{1,2,3\}$, and  $z=0$.
Inverse Laplace transformation yields the mean exit time as
\begin{equation}
\TT_b(x)=  \frac{2}{ p\lambda}+ \sum_{j=1}^3 \frac{ 1+q^{-1/3} e^{-2\pi ij/3}}{ 3 z_j} e^{z_j\varrho}.
 \end{equation}

Returning to the general case we see that
   the inversion of the
Laplace transform involves solving a  cubic
 equation, and though explicit formulas are available the resulting expression is
 awkward.  Still, the    large  $b$ limit
can be discerned with all generality.  To this end we note that by
appealing to Hurwitz's stability criteria it can be proven that
all three roots $z_j$, $j\in\{1,2,3\}$, of the denominator  in  expression
(\ref{FF}) |apart from $z=0$|  have negative real parts.  Hence,
evaluating the inversion integral by residues  we find
  \begin{eqnarray}
 \TT_b(x)=\frac{2}{ p\lambda}  + \sum_{j=1}^3C_j
 e^{z_j\varrho},
\end{eqnarray}
where $C_j$ are certain constants. Thus, letting $b\to\infty$ we
see that $\TT_b(x) \underset{b \to\infty}\to  \frac{2 }{ p\lambda} =\frac{\EE[\tau_n]}{1-\EE[\theta(J_n)]}$.

The evaluation of the correction to $T_b(x,r)$ does not present
particular difficulties. We leave it as an exercise to the
interested reader.

\section{Conclusions}

    We have analyzed the   mean exit time  for  a
general CTRW with drift. If the present  coincides with   a jump time we find that it  satisfies a
   certain integral equation whose solvability is analyzed. We consider next the generic  case when  the present is an
arbitrary instant and
the      history of the system is not available to the observer and only the present state is.
It turns out that the corresponding escape time can be obtained by incorporating an appropriate correction, which    can be described
in terms of the ``excess life'', a familiar object in renewal
theory.  We
   find that when the  drift and jump  components have  the same sign the equations that these objects satisfy  can be solved  in closed form
   via Laplace transformation, irrespective of the  distribution;
   otherwise, one must restrict to particular choices of the
sojourn-time distribution. The case   corresponding to the
classical Erlang distribution is analyzed in detail.    The more general case when jumps take both positive and negative signs is
also considered and solved under certain severe conditions.   We plan to generalize these ideas to a
more general class of waiting-time distributions  and   pinpoint
conditions that guarantee the  reducibility of the original formulation to simpler differential
equations.

The relevance  of  these results from a physical perspective is discussed in several
connections  of interest  including  possible
applications to risk, finance and distribution of energy in  optical systems, which will be the matter of future publications. We also point
out the relevance of the approach whenever   the  time between events  is
``large''
 or when the  arrival  times  are not physical
 observables.

\acknowledgments

The authors acknowledge support from MICINN under contracts No. 
FIS2008-01155-E, FIS2009-09689, and MTM2009-09676; from  Junta de Castilla y Le\'on, SA034A08; and Generalitat de Catalunya, 2009SGR417.

\appendix

\section{Some remarks on the continuum limit and its relationship with fractional diffusions}
In this appendix we sketch how the approach followed in this article, based in the use of renewal theory, compares with the most traditional one which relies on the previous computation of first-passage time PDFs.



In particular, we shall stress the connections of both techniques under the continuum limit approximation. This concept  |which  we will define properly in short| can be loosely identified with the limit in which both the mean sojourn time $\mu$ and the characteristic jump magnitude tend to zero |note that, by contrast, in this paper we have considered a situation wherein sojourn times are moderate or even large.

The major benefit of the continuum assumption is that it allows to obtain general results on the basis of limited knowledge of the jumping time and size PDFs, even when these distributions do not have all their moments well defined. The major drawback within our set-up is that as $\mu \to 0$ 
the variable $r$, the time elapsed from the last known jump, tends to zero as well, and the distinction between $T_b(x,r)$ and $\TT_b(x)$ becomes irrelevant. Therefore, any comparison between the two methods must be focused on how the object $\TT_b(x)$ is obtained.

Let us assume, for instance, that $\hat{\psi}(s)\sim 1-\mu s +o(s)$ and $\hat{h}(s)\sim 1-k |s|^{\alpha} +o(|s|^{\alpha})$ when $s\to 0$ for certain constants $\mu$, the mean sojourn time, and $k$. (Note however that it is also possible to consider the continuum limit in the case in which the mean sojourn time does not exist. We are just giving an illustrative example. For a more exhaustive analysis of this topic see~\cite{Metzler}.) 

To be more explicit, let us consider the case
\begin{equation}
h(u)=\frac{k}{2 \sqrt{\pi u^3}} e^{-\frac{k^2}{4 u}}\quad (u>0),
\end{equation}
so that,
\begin{equation}
\hat h(s)=e^{-k \sqrt{s}}\sim 1-k \sqrt{s} +o(\sqrt{s}).
\end{equation}
It is well known that for the CTRW process $S_t$, Eq.~(\ref{process}), the 
propagator 
$p(u,t)du\equiv\Pr\{u<S_t\leqslant u+du\}$ reads in the Laplace-Laplace domain
\begin{equation}
\hat p(s_1,s_2)=  \frac{1}{s_2}\frac{1-\hat{\psi}(s_2)}{1-\hat{\psi}(s_2)\hat{h}(s_1)}\sim\frac{\mu}{\mu s_2+k \sqrt{s_1}},
\end{equation}
as $s_{1,2}\to 0$. The continuum limit is recovered in this case by letting $\mu \to 0$, $k \to 0$ with $k/\mu \to K$ finite. The process arising after this limit fulfills the following fractional diffusion equation
\begin{equation}
\partial_t p(u,t) + K\ { }^{}_0 D_u^{\frac{1}{2}} p(u,t)=0, 
\end{equation}
where ${ }^{}_0 D_u^{\frac{1}{2}}$ is the Riemann-Liouville fractional operator of order $\frac{1}{2}$, and whose solution reads
\begin{equation}
p(u,t)=\frac{K t}{2 \sqrt{\pi u^3}} e^{-\frac{K^2 t^2}{4 u}}.
\end{equation}
Let us now define $\Pi_b(x,t)$ as the probability that the process $X_t$ has never touched $b$ when it started at $x< b$ at the initial time, i.e.  
\begin{equation}
\Pi_b(x,t)\equiv\Pr\{X_{t'}< b, 0\leqslant t' \leqslant t |X_0=x\}.
\end{equation}
In the present case, as $S_t$ is an increasing positive process and $v>0$, it can be computed by direct integration of $p(u,t)$: 
\begin{equation}
\Pi_b(x,t)=\int_0^{b-x-vt} p(u,t) du=\text{Erfc}\left(\frac{K^2 t^2}{2\sqrt{b-x-vt}}\right).
\end{equation}
Now one can obtain $\TT_b(x)$ through
\begin{eqnarray}
&&\TT_b(x)=\int_{0}^{\frac{b-x}{v}} t \partial_t\left[1- \Pi_b(x,t)\right] dt=\int_{0}^{\frac{b-x}{v}} \Pi_b(x,t) dt\nonumber\\
&&=\frac{4v}{\sqrt{\pi}K^2}\int_0^{\infty}\left(\xi \sqrt{\xi^2+\frac{K^2(b-x)}{v^2}}-\xi^2\right) e^{-\xi^2}d\xi \nonumber \\
&&=\frac{2}{K}\sqrt{\frac{b-x}{\pi}}+ \frac{v}{K^2} \left[ e^{\frac{K^2(b-x)}{v^2}}\text{Erfc}\left(\frac{K}{v}\sqrt{b-x}\right)-1\right]. \nonumber \\ \label{Ex_TTb}
\end{eqnarray}

Alternatively, if we use the direct approach followed in this paper our results imply that the mean arrival time follows by inversion of the Laplace expression in
Eq.~(\ref{FF_sol}) 
when $\mu \to 0$, $k \to 0$ with $k/\mu \to K$ finite
\begin{equation}
\hat{\FF}(s)=\frac{1}{v s^2}\frac{1-\hat{\psi}(s v)}{1-\hat{\psi}(s v)\hat{h}(s)} \rightarrow \frac{1}{vs^2 +Ks^{\frac{3}{2}}},
\end{equation}
that is \begin{equation}
\TT_b(x)=\frac{1}{2\pi i}
\int_{c-i\infty}^{c+i\infty} \frac{e^{ s(b-x)} }{vs^2 +Ks^{\frac{3}{2}}} d s.
\end{equation}
Upon evaluation of this integral the result (\ref{Ex_TTb}) is recovered.

\end{document}